\documentclass[11pt,a4paper]{article}
\usepackage{geometry}
\usepackage{amsmath,amsthm,amssymb,amsfonts}
\usepackage{amsfonts}
\newtheorem{remark}{Remark}
\usepackage{mathabx}
\usepackage[colorlinks=true,linkcolor=blue,citecolor=blue]{hyperref}

\usepackage{multirow}
\usepackage[boxed,linesnumbered]{algorithm2e}
\usepackage{algorithm2e}
\usepackage[table]{xcolor}
\usepackage{booktabs}
\usepackage{rotating, graphicx}
\usepackage{lscape}
\usepackage{longtable}
\usepackage{natbib}
\usepackage{appendix}
\usepackage{fixdif}
\usepackage{cleveref}
\crefname{equation}{}{}
\Crefname{algocf}{Algorithm}{Algorithms} 
\crefmultiformat{equation}{(#2#1#3)}{ and~(#2#1#3)}{, (#2#1#3)}{ and~(#2#1#3)}
\crefname{appsec}{Appendix}{Appendices}
\usepackage{todonotes}

\allowdisplaybreaks
\definecolor{LG}{gray}{0.9}
\definecolor{Gray}{gray}{0.85}

\newcommand{\mc}[1] {\mathcal{#1}}
\newtheorem{proposition}{Proposition}

\title{A Solution Method for Single-Source Capacitated Facility Location Problems Under Decision-Dependent Demand Uncertainty}
\author{Giovanni Pantuso \\
Department of Mathematical Sciences \\
University of Copenhagen \\
Universitetsparken 5, 2100, Copenhagen, Denmark\\
\href{mailto:gp@math.ku.dk}{gp@math.ku.dk}\\
Miguel Lejeune\\
The George Washington University\\
Washington, United States \\
\href{mailto:mlejeune@gwu.edu} {mlejeune@gwu.edu}\\
}
\date{}

\begin{document}
\maketitle
\begin{abstract}
  We present a model and an exact solution method for \textit{single-source capacitated facility location problems} with decision-dependent demand uncertainty.
  The model allows the joint probability distribution of customer demands to depend on the facility to which demand points are allocated.
  The number of possible decision-dependent distributions is therefore equal to the number of feasible allocations of demand points to (groups of) facilities.
  To address the exponential size of the model, we apply and extend a recent version of the L-Shaped method for stochastic programs with decision-dependent uncertainty. 
  In particular, we compare the original optimality cut with a new optimality cut based on a bilinear reformulation.
  We prove that, separated together with the original one, the new cut yields a tighter restriction, at the expense of introducing bilinear terms.
  \color{black} In addition, we prove an effective decision-independent valid inequality that significantly improves convergence.
  Extensive numerical tests on $180$ instances \color{black}demonstrate that the method, supported by the valid inequality, can find optimal or near-optimal solutions to problem instances with up to $6^{20}$ potential probability distributions.
  The tests also show that, while the new cut closes slightly more instances, the original cut remains preferable in terms of solution time and optimality gap.
  \color{black}
  
  {\bf Keywords:} Stochastic Programming, Endogenous Uncertainty, Facility Location, Decision-Dependent Uncertainty.
\end{abstract}

\section{Introduction}

Facility location decisions are significantly influenced by uncertainty \citep{LiuSWM22,CorS19,LapNS19,LapL94,Lou86,Sny06}.
The optimization literature is rich of methodology to address the resulting problems. The corresponding literature is thoroughly investigated and classified in the survey paper by \cite{Sny06} and recent books \citep{CorS19,SalW24}.

In several real-life applications, location decision problems are further complicated by additional constraints. 
A prime example, which is investigated in this paper, is the \textit{single-source} requirement. Upon making location decisions, the decision maker simultaneously assigns each customer to at most one capacitated facility and build the corresponding infrastructure. The allocated facility supplies its entire demand \citep{guastaroba2014heuristic}.
This requirement is rooted in several industrial applications.
For instance, in healthcare planning, decision makers locate clinics and simultaneously assign population zones to a single clinic.
Similarly, when locating retail stores, firms may assign geographical zones to a single retail, which covers its entire demand.
In the telecommunications industry, service providers simultaneously decide where to place optical fiber aggregation nodes and which households to connect to each installed node. Each household is connected to exactly one node, requiring fixed investments in trenching, cabling, and installation. In water supply system planning, decision makers determine which treatment plants to build and assign each demand zone to a single plant. This assignment entails constructing pipelines and pressure-control infrastructure.
The single-source constraint makes the problem significantly harder to solve, see \cite{klose2005facility}. For example, assume that a set of open facilities has been chosen. If multiple sources are allowed, the decision problem of optimally supplying customers from those facilities forms a transportation problem. However, if single-sourcing is required, the corresponding assignment problem forms a special \textit{Generalized Assignment Problem} (GAP) which is NP-hard, see \cite{fisher1986multiplier}.

In some applications, demand uncertainty may depend on the facility serving a given area. In healthcare planning, demand for healthcare services may vary with the accessibility and quality of the assigned clinic. Zones served by nearby or better-equipped clinics may exhibit higher utilization rates.
In retail network design, the demand generated by a zone may depend on the attractiveness of the assigned store. Factors, such as distance, product assortment, and service quality, can all influence customer demand and therefore affect both location and assignment decisions.
In the planning of emergency and public services, such as ambulance, vaccination, or social-service systems, utilization may increases when services become easier to access. Assigning a region to a closer facility may induce additional demand.
Finally, in the planning of transportation systems, travelers may use a transit hub, station, or shared-mobility facility more often when assigned to a convenient or high-capacity location.

Despite the abundance of methodological approaches, the likely influence of location decisions on the uncertainties has not been much 
\color{black}
investigated. 
To our knowledge, the only available studies, which address the impact of location decisions on the uncertainty are those of \cite{BasAS21} and \cite{Pan25}.
\cite{BasAS21} study a facility location problem where demand is endogenously influenced by location decisions.
Particularly, the mean and variance of the demand distribution is a function of facility-location decisions.
The authors propose a distributionally-robust optimization model to address misspecifications in the demand distribution.
That is, location decisions determine an ambiguity set, which contains the potential demand distributions that match the first two moments.
The objective of the model is then to minimize the facility-location costs and the worst-case expected costs 
of transportation and unsatisfied demand.
The authors show that the model can be reformulated as a 
mixed-integer linear programming 
\color{black}
(MILP) problem in the special case where demand mean and variance are piecewise linear functions of location decisions.
In an extensive computational study, the authors demonstrate the superiority of the proposed model to traditional robust optimization and stochastic programs.
\cite{Pan25} studies a facility location problem where all parameters of the problem may be affected by uncertainty. The author presents a general mathematical model in the form of a stochastic program with decision-dependent uncertainty. The central assumption is that the distribution of the data depends on the facilities (or groups of facilities) opened in the first-stage solution. A consequence of this is that there exist only finitely many possible distributions that may be induced, one for each feasible facility location decision. However, 
 the number of feasible solutions -- and corresponding distributions -- is $\mathcal{O}(2^I)$ with $I$ denoting the number of facilities.
The author shows that the problem can be formulated as an MILP problem with exponentially many constraints. To solve the problem, the author proposed an exact decomposition method. Exact solutions are found for instances with up to $2^{17}$ distributions. 

In the case of decision-dependent uncertainty, the distribution of demand can no longer be meaningfully assumed to depend on the choice of open facilities, as done in \cite{Pan25}. 

When each customer’s demand is satisfied entirely by a single source, it is more realistic to assume that the demand distribution depends on the facility to which the customer is assigned. Consequently, if $I$ denotes the number of facilities and $J$ the number of customers, the number of allocation decisions
-- and corresponding distributions -- is $\mathcal{O}(I+1)^J$, see \Cref{sec:model}. This number is, in general, much larger than $2^I$ as in the case of \cite{Pan25}. 
We address this case in this article.

At a more general level, stochastic programs with decision-dependent uncertainty are notoriously challenging.
While the research literature offers several practical applications, see e.g., \cite{VisSF04,PeeSGV10,LJ1,TonFR12,EscGMU18,LJ2,LJ3}, general-purpose methodology remains sparse. 
In particular, \cite{JonWW98} study decision problems where decisions have an impact both on the probability of future random events, as well as on timing of their resolution.  
The authors assume that the set of probability measures which can be enforced by decisions is finite, as is the case in our paper. The authors show that, for moderately-sized instances, the problem can be solved by enumerating the possible distributions and solving the corresponding stochastic program.
\cite{Pan21} extends the model of \cite{JonWW98} by allowing decisions at all stages to determine the probability measure. The author proposes a formulation that avoids explicit non-anticipativity constraints, see e.g., \cite{ApaG17,HooM16,MooM18}.
\cite{HelBT18} discuss several ways 

to model the interplay between decisions and uncertainties.
\color{black}
In their work, it is evident how the resulting formulations are typically nonlinear. Their characteristics heavily depend on the nature of the relationship between decisions and distributions, severely complicating the development of general-purpose solution approaches. 
Finally, \citet{PanH25} propose an extension of the L-Shaped method for two-stage stochastic programs with decision-dependent uncertainty. The method is based on identifying the enforced distribution for each candidate first-stage solution prior to separate decision-dependent cuts.
It relies on the assumption that there exists finitely many probability distributions that may be enforced.
This assumption naturally holds for many combinatorial optimization problems, such as the one studied in this paper.


In this paper we consider and extend the \textit{single-source capacitated stochastic facility location problem} presented by \cite{LapL94} (see also \cite{HinK99,AhuOPSS04,Lin09,YanCC12,LapNS19} for variants of the problem).
In this problem, each demand point may be assigned to at most one facility.
Hence, first-stage decisions determine both the open facilities and the assignment of demand points to facilities.
A fixed cost is paid upon opening a facility and creating the link between a demand point and a facility.
Demand may be satisfied only through open links.
Each point has a stochastic demand, which is revealed only after facilities and links have been open.
For this reason, the demand does not need to be met exactly. Rather, a revenue per unit of satisfied demand is introduced.
The setup of \cite{LapL94} is based on the assumption that 
the demand distribution is known and independent of decisions.
We relax this (restrictive) assumption 
by allowing the probability distribution of the demands to depend on the facility to which the demand point is allocated.
This, in turn, may allow modeling situations where demand is influenced by e.g., the distance to the provider, the cost born for each unit of demand satisfied, or characteristics of the providers. 
This generates a problem where there exists one probability distribution for each possible allocation of demand points to facilities. As we will see, this number is exponential in the number of customers.


The contribution of this paper is as follows:
\begin{itemize}
\item We propose a stochastic programming model for a \textit{single-source capacitated facility location problem} under endogenous stochastic demand.
  The model allows the distribution of demand to depend on the allocation of demand points to facilities.
  The proposed model is general in the sense that it models a range of scenarios having the original problem presented by \cite{LapL94} at one extreme and, at the other extreme,
  a problem where there are as many distributions as there are allocations of demand points to facilities.
\item To solve the problem, we present an exact cutting-plane algorithm\color{black}, which can be regarded as a version of the L-Shaped method for stochastic programs with decisions-dependent uncertainty developed by \cite{PanH25}.
  In particular, we develop a new general optimality cut and prove that, separated together with the cut of \cite{PanH25}, it yields a tighter restriction.
  The tightening requires the addition of bilinear terms, which can be linearized exactly.
  This, however, introduces a number of linearization variables that grows quadratically with the number of demand points. \color{black}
  In addition, we prove a very tight, problem-specific, valid inequality, which significantly increases the performance of the method.
\item We present results of extensive numerical tests.
  These illustrate that the method can solve instances with up to $6^{20}$ potential probability distributions.
  They also show that the new cut closes slightly more instances than the original one, but the original cut remains preferable overall, both in terms of solution time and optimality gap attained.
  We explain that the bilinear terms make each round of separation more expensive leading to the generation of fewer cuts within the same time limit.
  This illustrates that a cut that is provably tighter is not necessarily more effective once the cost of separating it is taken into account.\color{black}
\end{itemize}
In the remainder of this article we begin by presenting the model in \Cref{sec:model}.
Following, we describe the solution method in \Cref{sec:method}.
In \Cref{sec:experiments}, we present the results of numerical experiments. Finally, we draw conclusions in \Cref{sec:conclusions}.

\section{Model}\label{sec:model}
Let $\mc{I}$ be the set of potential facility locations and $\mc{J}$ be the set of demand points to serve.
Opening a facility at location $i$ determines a fixed cost $F_i$. Opening a connection between facility $i$ and demand point $j$ determines a fixed cost $C_{ij}$.
Each demand point has uncertain demand. We model the demand as a $\vert\mc{J}\vert$-dimensional real-valued random variable $\xi:(\Omega,\mc{F})\to(\mathbb{R}^{|\mc{J}|},\mc{B}^{|\mc{J}|})$ on a measurable space $(\Omega,\mc{F})$, where $\mc{B}^{|\mc{J}|}$ is the Borel $\sigma$-algebra on $\mathbb{R}^{|\mc{J}|}$. 
We allow the probability distribution of $\xi$ on $\mathbb{R}^{\vert\mc{J}\vert}$ to depend on the location of the facility to which each demand point has been assigned.

Let $x:=(x_i)_{i\in\mc{I}}\in\{0,1\}^{|\mc{I}|}$ denote the vector of facility location decisions and let $y:=(y_{ij})_{i\in\mc{I},j\in\mc{J}}\in\{0,1\}^{|\mc{I}|\times |\mc{J}|}$ represent the assignments of demand points to facilities.
To allow for a rather general setup, we assume that there exists a partition $(\mc{I}_z)_{z=1,\ldots,Z}\subset 2^{\mc{I}}$ of $\mc{I}$, into $Z$ (e.g., geographical) subsets, where $2^{\mc{I}}$ is the power set of $\mc{I}$ with $Z\leq|\mc{I}|$ obviously.
The binary variable vector $u:=(u_{jz})_{j\in\mc{J},z=1,\ldots,Z+1}\in\{0,1\}^{\vert\mc{J}\vert (Z+1)}$ represents the allocations of demand points to subsets, i.e., $u_{jz}$ implies that demand point $j$ has been assigned to a facility in subset $z$.
Note that we introduce a dummy subset $\mc{I}_{Z+1}:=\emptyset$ \color{black} to which we allocate demand points that have not been assigned to any facility. Hence, the subsets $(\mc{I}_z)_{z=1,\ldots,Z+1}$ \color{black} form likewise a partition of $\mc{I}$. 

We let the assignment of demand points to subsets determine the joint distribution of $\xi$ on $\mathbb{R}^{\vert\mc{J}\vert}$.
This gives $(Z+1)^{|\mc{J}|}$ potential distributions, one for each possible allocation of demand points to subsets.
Observe that this allows us to model a vast array of potential configurations.
At one end, we find the original case where all locations are in the same subset and the distribution is not affected by allocation decisions.
At the other end, we find the case where $Z=(|\mc{I}|+1)$ and the joint demand distribution depends on the individual facilities to which demand points have been allocated.
By acting on the number of subsets $Z$, one obtains that the distribution depends on the subset of facilities to which demand points are assigned.
To gain some practical intuition, one may think of subsets as, e.g., geographical zones or facility types.
Let then $\mc{D}$ be the set of all possible joint distributions of demand that may materialize: $\vert\mc{D}\vert=(Z+1)^{\vert\mc{J}\vert}$. This gives a probability distribution $\mathsf{P}_d$ for each possible allocation decision $d\in\mc{D}$.
One may for example think of a distribution parametrized by a function of $u$, i.e., $\mathsf{P}^\theta$, with $\theta:=f(u)$.
Given that $u$ can take only finitely many values, one can enumerate the possible specifications of $\mathsf{P}^\theta$ as $\mathsf{P}_1,\ldots,\mathsf{P}_{\vert\mc{D}\vert}$.
Let each component $\delta_d$ of the binary variable vector $\delta:=(\delta_d)_{d\in\mc{D}}\in\{0,1\}^{\vert\mc{D}\vert}$ take value one if distribution $d$ is enforced by the allocation decisions, and be equal to $0$ otherwise. 
Finally, we define a map $\zeta:\mc{J}\times \mc{D}\to \{1,\ldots,Z+1\}$, which returns the zone to which a demand point $j$ must be assigned in order for the distribution $d$ to accurately describe the uncertain demand.
Having this notation in place, the \textit{single-source facility location problem with decision-dependent stochastic demand} can be expressed as follows:
\begin{subequations}
  \label{eq:flp}
  \begin{align}
    \label{eq:flp:1st:obj}\min &\sum_{i\in\mc{I}}F_ix_i+\sum_{i\in\mc{I}}\sum_{j\in\mc{J}}C_{ij}y_{ij}-\sum_{d\in\mc{D}}\delta_d\mathbb{E}_{\mathsf{P}_{d}}Q(y,\xi)\\
    \label{eq:flp:1st:1}\text{s.t.~}&y_{ij}\leq x_i &\forall i\in\mc{I},j\in\mc{J}\\
    \label{eq:flp:1st:2} &\sum_{i\in\mc{I}}y_{ij}\leq 1 &\forall j\in\mc{J}\\
    \label{eq:flp:1st:assignment1}           &\sum_{i\in\mc{I}_z}y_{ij}=u_{jz}&\forall j\in\mc{J},z=1,\ldots,Z\\
    \label{eq:flp:1st:assignment2}            &\sum_{z=1}^{Z+1}u_{jz} = 1 &\forall j\in\mc{J}\\
    \label{eq:flp:1st:impl1}                               &\sum_{j\in\mc{J}}u_{j,\zeta(j,d)}\leq \vert \mc{J}\vert +\delta_d -1&\forall d\in\mc{D}\\
    \label{eq:flp:1st:impl2}                               &\sum_{j\in\mc{J}}u_{j,\zeta(j,d)}\geq \vert \mc{J}\vert \delta_d &\forall d\in\mc{D}\\
    \label{eq:flp:1st:3}&x_{i}\in\{0,1\} &\forall i\in\mc{I}\\
    \label{eq:flp:1st:4}&y_{ij}\in\{0,1\} &\forall i\in\mc{I},j\in\mc{J}\\
    \label{eq:flp:1st:5}&u_{jz}\in\{0,1\} &\forall j\in\mc{J},z=1,\ldots,Z+1\\
    \label{eq:flp:1st:6}&\delta_{d}\in\{0,1\} &\forall d\in\mc{D}
  \end{align}
\end{subequations}
The objective function sums the costs of opening facilities and links to demand points, and accounts for the stochastic revenue made as a consequence of such decisions.
Here $\mathbb{E}_{\mathsf{P}_{d}} Q(y,\xi) =\int Q(y,\xi)\d \mathsf{P}_d$ represents the expected second-stage revenue according to distribution $\mathsf{P}_d$. We assume this expectation to be finite for all $y$ and $d\in\mc{D}$.
The second-stage function $Q(y,\xi)$ will be introduced later.

Constraints \eqref{eq:flp:1st:1} state that a connection can be open only if the corresponding facility is open.
Constraints \eqref{eq:flp:1st:2} state that each demand point can be connected to at most one facility.
Constraints \eqref{eq:flp:1st:assignment1} set $u_{jz}=1$ whenever demand point $j$ has been assigned to a facility in set $z$. Constraints \eqref{eq:flp:1st:assignment2} ensure that
each demand point is either assigned to a real subset or to the dummy set $Z+1$. Recall that the assignment of demand points to subsets determines the demand distribution.
Here constraints \eqref{eq:flp:1st:impl1} model the implication $\sum_{j\in\mc{J}}u_{j,\zeta(j,d)}=\vert\mc{J}\vert \implies \delta_d =1$ while \eqref{eq:flp:1st:impl2} implement
$\delta_d=1 \implies \sum_{j\in\mc{J}}u_{j,\zeta(j,d)}=\vert\mc{J}\vert$. This double implication ensures that distribution $\mathsf{P}_d$ is enforced if and only if the corresponding allocation of demand points to subsets, according to $\zeta$, holds.

Finally, we represent the second-stage problem as follows:
\begin{subequations}
  \begin{align}
    \label{eq:flp:2st:obj}Q(y,\xi) = & \max \sum_{i\in\mc{I}}\sum_{j\in\mc{J}}R_{ij}w_{ij}\\
    \label{eq:flp:2st:1}\text{s.t.~}& w_{ij}\leq \xi_j y_{ij}&\forall i\in\mc{I},j\in\mc{J}\\
    \label{eq:flp:2st:2}        &\sum_{j\in\mc{J}}w_{ij}\leq M_i&\forall i\in\mc{I}\\
         \label{eq:flp:2st:3}        &w_{ij}\geq 0&\forall i\in\mc{I},j\in\mc{J}
  \end{align}
\end{subequations}
Here $w_{ij}$ represents the demand of demand point $j$ satisfied by facility $i$.
The objective is to maximize revenues, where $R_{ij}\geq 0$ is the per-unit revenue for the demand of demand point $j$ satisfied from facility $i$.
The constraints ensure that demands $\xi_j$ and facility capacities $M_i\geq 0$ are not exceeded.

Observe that model \eqref{eq:flp} is in general nonlinear due to the product in the last term of the objective function.
For discrete probability distributions defined with finitely many scenarios, we obtain a bilinear optimization problem due to the product of variables $\delta_d$ and $w_{ij}$ in the objective function.
These bilinear terms can be linearized exactly using the McCormick method \citep{McC76}, 
which yields an MILP problem\color{black}. 
Nevertheless, the number of constraints \eqref{eq:flp:1st:impl1}-\eqref{eq:flp:1st:impl2} is exponential in $\vert\mc{J}\vert$.

\section{Solution Algorithm}\label{sec:method}
We develop a cutting-plane algorithm that instantiates and extends the L-Shaped method for stochastic programs with decision-dependent uncertainty \citep{PanH25}. In addition, we provide a novel, generally-applicable, optimality cut, which is tighter than one developed by \cite{PanH25}, but requires dealing with bilinear terms that can be linearized exactly. 

In \Cref{sec:ls}, we first present the L-Shaped method for stochastic programs with decision-dependent uncertainty. Then, in \Cref{sec:bilinear_cuts}, we introduce the new optimality cuts. In \Cref{sec:validinequality}, we introduce a problem-specific valid inequality.

\subsection{The L-Shaped Method with Decision-Dependent Distributions}\label{sec:ls}
In this section, we extend the method of \cite{PanH25} to our problem. 
We assume that $\xi$ is a discrete random variable so that distributions $\mathsf{P}_d$, $d\in\mc{D}$, are step functions. 
Let $\xi_{ds}$ and $\pi_{ds}$ respectively denote the realization of the random variable $\xi$ and its probability under scenario $s$ of distribution $\mathsf{P}_d$ and let $\mc{S}_d$ be the distribution-specific set of scenarios characterizing~$\mathsf{P}_d$.

We first obtain the epigraph reformulation of problem \eqref{eq:flp} by adding a decision variable~$\mu$:
\begin{subequations}
  \label{eq:flpref}
  \begin{align}
    \label{eq:flpref:1st:obj}\min  &\sum_{i\in\mc{I}}F_ix_i+\sum_{i\in\mc{I}}\sum_{j\in\mc{J}}C_{ij}y_{ij}-\mu\\
\text{s.t.~}&\cref{eq:flp:1st:1,eq:flp:1st:2,eq:flp:1st:assignment1,eq:flp:1st:assignment2,eq:flp:1st:impl1,eq:flp:1st:impl2,eq:flp:1st:3,eq:flp:1st:4,eq:flp:1st:5,eq:flp:1st:6}\\
    \label{eq:flpref:1st:Qx}&\mu\leq\sum_{d\in\mc{D}}\delta_d\mathbb{E}_{\mathsf{P}_{d}}Q(y,\xi)\\
    &\mu\geq 0             
  \end{align}
\end{subequations}
Observe that, given that the subsets $(\mc{I}_z)_{z=1,\ldots,Z+1}$ \color{black} form a partition of $\mc{I}$, only one $\delta_d$ can take value $1$.
Furthermore, observe that $Q(y,\xi)$ is always bounded from above under the reasonable assumption that $M_i$ and $\xi_j$ are finite.
In addition, since $(0)_{i\in\mc{I},j\in\mc{J}}$ is always feasible, it follows that $Q(y,\xi)\geq 0$ for all $y$ and $\xi$ and that there exists an upper bound $U$ on the second-stage expected revenue given by: 
$$U\geq \max_{d\in\mc{D}}\max_{y\in\{0,1\}^{\vert\mc{I}\vert\vert\mc{J}\vert}}\mathbb{E}_{\mathsf{P}_d}Q(y,\xi) \ . $$

\Cref{prop:ub} provides a closed-form expression for $U$.
\begin{proposition}\label{prop:ub}
  Let
  $$U_1:=\sum_{j\in\mc{J}}\max_{i\in\mc{I}}\left\{R_{ij}\min\left\{M_i,\max_{d\in\mc{D}}\mathbb{E}_{\mathsf{P}_d}[\xi_j]\right\}\right\}$$
  and
  $$U_2:=\left(\max_{i\in\mc{I},j\in\mc{J}}R_{ij}\right)\sum_{i\in\mc{I}}M_i$$
  Then $U:=\min\{U_1,U_2\}$ satisfies $U\geq\mathbb{E}_{\mathsf{P}_d}Q(y,\xi)$ for every $d\in\mc{D}$ and every $y$ satisfying \eqref{eq:flp:1st:2} and \eqref{eq:flp:1st:4}.
\end{proposition}
\begin{proof}
  It suffices to prove the bound for $U_1$ and for $U_2$ separately.

  Consider $U_2$.
  Let $w$ be feasible for $Q(y,\xi)$.
  Since $w\geq 0$ and $R_{ij}\leq\max_{i,j}R_{ij}$ for all $i,j$, constraints \eqref{eq:flp:2st:2} give
  $$\sum_{i\in\mc{I}}\sum_{j\in\mc{J}}R_{ij}w_{ij}\leq\left(\max_{i,j}R_{ij}\right)\sum_{i\in\mc{I}}\sum_{j\in\mc{J}}w_{ij}\leq\left(\max_{i,j}R_{ij}\right)\sum_{i\in\mc{I}}M_i=U_2 $$
  Hence $Q(y,\xi)\leq U_2$ for every $y$ and every realization $\xi$, and therefore $\mathbb{E}_{\mathsf{P}_d}Q(y,\xi)\leq U_2$ for every $d$.

  Consider now $U_1$. Fix $d\in\mc{D}$ and $y$.
  The function $Q(y,\cdot)$ is the optimal value of a maximization linear program in which $\xi$ enters only the right-hand side of \eqref{eq:flp:2st:1}.
  It follows that $Q(y,\cdot)$ is concave in $\xi$. Jensen's inequality gives
  $$\mathbb{E}_{\mathsf{P}_d}Q(y,\xi)\leq Q\left(y,\mathbb{E}_{\mathsf{P}_d}[\xi]\right)$$
  Let $w$ be feasible for $Q(y,\mathbb{E}_{\mathsf{P}_d}[\xi])$.
  By \eqref{eq:flp:2st:1}, $w_{ij}\leq\mathbb{E}_{\mathsf{P}_d}[\xi_j]y_{ij}$, and by \eqref{eq:flp:2st:2} together with $w\geq 0$, $w_{ij}\leq\sum_{j'\in\mc{J}}w_{ij'}\leq M_i$, so that, for all $i\in\mc{I}$ and $j\in\mc{J}$ we have
  $$w_{ij}\leq\min\left\{M_i,\mathbb{E}_{\mathsf{P}_d}[\xi_j]\right\}y_{ij}$$
  Therefore
  $$Q\left(y,\mathbb{E}_{\mathsf{P}_d}[\xi]\right)\leq\sum_{j\in\mc{J}}\sum_{i\in\mc{I}}R_{ij}\min\left\{M_i,\mathbb{E}_{\mathsf{P}_d}[\xi_j]\right\}y_{ij}$$
  By \eqref{eq:flp:1st:2} and \eqref{eq:flp:1st:4}, for each $j\in\mc{J}$ at most one $y_{ij}$ equals one.
  Hence, the inner sum is bounded by its largest term, that is
  $$\sum_{i\in\mc{I}}R_{ij}\min\left\{M_i,\mathbb{E}_{\mathsf{P}_d}[\xi_j]\right\}y_{ij}\leq\max_{i\in\mc{I}}\left\{R_{ij}\min\left\{M_i,\mathbb{E}_{\mathsf{P}_d}[\xi_j]\right\}\right\}$$
  Finally, $t\mapsto\min\{M_i,t\}$ is non-decreasing, so $\min\{M_i,\mathbb{E}_{\mathsf{P}_d}[\xi_j]\}\leq\min\{M_i,\max_{d'\in\mc{D}}\mathbb{E}_{\mathsf{P}_{d'}}[\xi_j]\}$ for every $d$.
  Summing over $j\in\mc{J}$ yields $\mathbb{E}_{\mathsf{P}_d}Q(y,\xi)\leq U_1$.
  Since $d$ and $y$ were arbitrary the claim follows.
\end{proof}
\color{black}
As a result, we can rewrite \eqref{eq:flpref:1st:Qx} as follows:
\begin{align}
  \label{eq:flpref:1st:Qxb}& \mu\leq\mathbb{E}_{\mathsf{P}_{d}}Q(y,\xi)+U(1-\delta_d)& \forall d\in\mc{D}
\end{align}
Note that the variable $\delta_d$, which takes 
value $1$ when distribution effectively $d$ is enforced,
will set the tightest bound on $\mu$ and make constraints \eqref{eq:flpref:1st:Qxb} equivalent to \eqref{eq:flpref:1st:Qx}.

Since $\xi$ is discrete, we may then write \eqref{eq:flpref:1st:Qxb} as follows:
\begin{align}
  \label{eq:flpref:1st:Qxc}& \mu\leq\sum_{s\in\mc{S}_d}\pi_{ds}Q(y,\xi_{ds})+U(1-\delta_d)& \forall d\in\mc{D} \ .
\end{align}
Since each $Q(y,\xi_{ds})$ is always feasible and bounded, strong duality implies that the optimal solution to the dual problem is always attained at one of the extreme points of the set:
$$\left\{(\rho,\gamma)\in\mathbb{R}^{\vert\mc{I}\vert \vert\mc{J}\vert }_{\geq 0}\times \mathbb{R}^{\vert\mc{I}\vert}_{\geq 0}\left\vert \gamma_i+\rho_{ij}\geq R_{ij}, i\in\mc{I},j\in\mc{J}\right.\right\} \ .$$
Let $(\rho^k,\gamma^k)$ be the $k$-th extreme point and $\mc{R}=\{(\rho^1,\gamma^1),\ldots, (\rho^{K},\gamma^{K})\}$ be the collection of such extreme points.
Hence, for each $d\in\mc{D}$ and $s\in\mc{S}_{d}$, we have:
$$Q(y,\xi_{ds})=\min_{(\rho,\gamma)\in\mc{R}}\left\{\sum_{i\in\mc{I}}\sum_{j\in\mc{J}}\rho_{ij}\xi_{jds}y_{ij}+\sum_{i\in\mc{I}}\gamma_iM_i\right\} \ ,$$
which allows the rewriting of \eqref{eq:flpref:1st:Qxc} as
\begin{align}
  \label{eq:flpref:1st:Qxd}& \mu\leq\sum_{s\in\mc{S}_d}\pi_{ds}\left(\sum_{i\in\mc{I}}\sum_{j\in\mc{J}}\rho_{ij}^{sd}\xi_{jds}y_{ij}+\sum_{i\in\mc{I}}\gamma_i^{sd}M_i\right)+U(1-\delta_d), \ \forall d\in\mc{D}, (\rho^{sd},\gamma^{sd})_{s\in\mc{S}_d}\in \times_{s\in\mc{S}_d}\mc{R},
\end{align}
where $\times_{s\in\mc{S}_d}$ denotes the cartesian product.
Note that the last term in \eqref{eq:flpref:1st:Qxd} models the implication:
\begin{align*}
  \delta_d=1\implies \mu\leq\sum_{s\in\mc{S}_d}\pi_{ds}\left(\sum_{i\in\mc{I}}\sum_{j\in\mc{J}}\rho_{ij}^{sd}\xi_{jds}y_{ij}+\sum_{i\in\mc{I}}\gamma_i^{sd}M_i\right) & \forall (\rho^{sd},\gamma^{sd})_{s\in\mc{S}_d}\in \times_{s\in\mc{S}_d}\mc{R} \ .
\end{align*}
As presented next, the same implication can be modeled using only the variables $u$.
Recall the map $\zeta:\mc{J}\times \mc{D}\to \{1,\ldots,Z+1\}$, which returns the subset $z$ to which $j$ has to be assigned in order for distribution $\mathsf{P}_d$ to apply.
Recalling that $\vert\mc{J}\vert-\sum_{j\in\mc{J}}u_{j,\zeta(j,d)}=0\implies \sum_{j\in\mc{J}}u_{j,\zeta(j,d)}=\vert\mc{J}\vert\implies \delta_d=1$, we can rewrite  \eqref{eq:flpref:1st:Qxd} as
\begin{align}
  \label{eq:flpref:1st:Qxe}
  &\mu\leq\sum_{s\in\mc{S}_d}\pi_{ds}\left(\sum_{i\in\mc{I}}\sum_{j\in\mc{J}} \rho_{ij}^{sd}\xi_{jds}y_{ij}+\sum_{i\in\mc{I}}\gamma_i^{sd}M_i\right)+U\left(\vert\mc{J}\vert-\sum_{j\in\mc{J}}u_{j,\zeta(j,d)}\right), \ \forall d\in\mc{D}, (\rho^{sd},\gamma^{sd})_{s\in\mc{S}_d}\in \times_{s\in\mc{S}_d}\mc{R},
\end{align}
where we have:
\begin{align*}
  \vert\mc{J}\vert-\sum_{j\in\mc{J}}u_{j,\zeta(j,d)}=0\implies \mu\leq\sum_{s\in\mc{S}_d}\pi_{ds}\left(\sum_{i\in\mc{I}}\sum_{j\in\mc{J}}\rho_{ij}^{sd}\xi_{jds}y_{ij}+\sum_{i\in\mc{I}}\gamma_i^{sd}M_i\right) & \forall (\rho^{sd},\gamma^{sd})_{s\in\mc{S}_d}\in \times_{s\in\mc{S}_d}\mc{R}.
\end{align*}
This effectively makes the variables $\delta$ redundant. Hence we arrive at the following reformulation of \eqref{eq:flp}, which we will refer to as the Master Problem (MP):

\begin{align}
  \label{eq:mp}\min \left\{\sum_{i\in\mc{I}}F_ix_i+\sum_{i\in\mc{I}}\sum_{j\in\mc{J}}C_{ij}y_{ij}-\mu \vert  \eqref{eq:flp:1st:1},\eqref{eq:flp:1st:2}, \eqref{eq:flp:1st:assignment1}, \eqref{eq:flp:1st:assignment2}, \eqref{eq:flp:1st:3},   \eqref{eq:flp:1st:4},   \eqref{eq:flp:1st:5},\eqref{eq:flpref:1st:Qxe}\right\} \ . \tag{MP}
\end{align}

Observe that, compared to \eqref{eq:flp}, \eqref{eq:mp} does not include constraints \eqref{eq:flp:1st:impl1} and \eqref{eq:flp:1st:impl2} and the $\delta_d$ variables \eqref{eq:flp:1st:6}, which are specified for all elements in $\mc{D}$. On the other hand, it includes constraints \eqref{eq:flpref:1st:Qxe}, which require enumerating the extreme points of the dual of problem $Q(y,\xi_{ds})$ for each $d$.

We propose to solve \eqref{eq:mp} by means of a cutting-plane algorithm, which relaxes constraints \eqref{eq:flpref:1st:Qxe} 
and adds cuts when found violated. 
Hence, at each iteration we solve the following \textit{Relaxed Master Problem} (RMP):
\begin{align}
  \label{eq:rmp}\min \left\{\sum_{i\in\mc{I}}F_ix_i+\sum_{i\in\mc{I}}\sum_{j\in\mc{J}}C_{ij}y_{ij}-\mu \vert  \eqref{eq:flp:1st:1},\eqref{eq:flp:1st:2}, \eqref{eq:flp:1st:assignment1}, \eqref{eq:flp:1st:assignment2},  \eqref{eq:flp:1st:3},   \eqref{eq:flp:1st:4},   \eqref{eq:flp:1st:5}\right\} \ .\tag{RMP}
\end{align}
Observe that \eqref{eq:rmp}, as displayed, contains no constraint bounding $\mu$ from above.
This would result in an unbounded problem when no optimality cut is present.
We therefore impose $\mu\leq U$, which is valid by definition of $U$ and never binding at an optimal solution of \eqref{eq:mp}.\color{black}

The procedure is sketched in \Cref{alg:ls}. At every iteration $v$, \eqref{eq:rmp} is solved with the currently available constraints \eqref{eq:flpref:1st:Qxe} (initially none).
Given a solution $(x^v,y^v,u^v,\mu^v)$ for \eqref{eq:rmp}, we identify the distribution $\mathsf{P}_{d^v}$ 
enforced by $u^v$, that is by the allocation of demand points to subsets of $\mc{I}$.
We compute next the expected second-stage revenue associated to that probability distribution.
Observe at line 6 of the pseudo-code that
\begin{align*}
  \mathbb{E}_{\mathsf{P}_{d^v}}[Q(y^v,\xi)] \; &= \; \sum_{s\in\mc{S}_d}\pi_{d^vs}Q(y^v,\xi_{d^vs})
&=\sum_{s\in\mc{S}_{d^v}}\pi_{d^vs}\left(\sum_{i\in\mc{I}}\sum_{j\in\mc{J}}\rho_{ij}^{sd^vv}\xi_{jd^vs}y^v_{ij}+\sum_{i\in\mc{I}}\gamma_i^{sd^vv}M_i\right) \ ,
\end{align*}
where $\rho^{sd^vv}_{ij}$ and $\gamma_i^{sd^vv}$ are the dual solution to $Q(y^v,\xi_{d^vs})$.
Hence, computing the expectation amounts to solving the second-stage problems only for the scenarios of distribution $d^v$.
If $\mu^v$ holds the correct value of the expectation, the algorithm has converged to an optimal solution.
Otherwise, we generate a constraint \eqref{eq:flpref:1st:Qxe} for $d=d^v$.
Particularly, if $\mu^v>\mathbb{E}_{\mathsf{P}_{d^v}}[Q(y^v,\xi)]$, we identify a violated inequality \eqref{eq:flpref:1st:Qxe} as:
\begin{equation}\label{eq:oc}
\mu\leq \sum_{s\in\mc{S}_{d^v}}\pi_{d^vs}\left(\sum_{i\in\mc{I}}\sum_{j\in\mc{J}}\rho_{ij}^{sd^vv}\xi_{jd^vs}y_{ij}+\sum_{i\in\mc{I}}\gamma_i^{sd^vv}M_i\right)+U\left(\vert\mc{J}\vert-\sum_{j\in\mc{J}}u_{j,\zeta(j,d^v)}\right) \ . \tag{OC}
\end{equation}
Observe that $(x^v,y^v,u^v,\mu^v)$ violates \eqref{eq:oc} by effect of having $\vert\mc{J}\vert=\sum_{j\in\mc{J}}u^v_{j,\zeta(j,d^v)}$.
Hence, we add \eqref{eq:oc} to \eqref{eq:rmp} as an \textit{optimality cut} that renders $(x^v,y^v,u^v,\mu^v)$ infeasible in the subsequent iterations.
Observe also how \eqref{eq:oc} is redundant for distributions other than $d^v$.
\begin{proposition}\label{prop:convergence}
  \Cref{alg:ls} converges in a finite number of iterations.
  \begin{proof}
    As we have shown, \eqref{eq:mp} is a reformulation of \eqref{eq:flp}. Recovering all (finitely many) inequalities \eqref{eq:flpref:1st:Qxe} is a sufficient condition for convergence.
    Therefore, the algorithm terminates in a finite number of iterations.
  \end{proof}
\end{proposition}

\begin{algorithm}
  \caption{Decision-Dependent L-Shaped Method for \eqref{eq:flp}.}\label{alg:ls}
  \texttt{solved}$\gets$\texttt{false}\;
  $v\gets 0$\;
  \While{not \texttt{solved}}{
    Solve \eqref{eq:rmp}. Let $(x^v,y^v,u^v,\mu^v)$ be its optimal solution\;
    Identify the distribution $\mathsf{P}_{d^v}$, $d^v\in\mc{D}$, enforced by allocation decisions $u^v$\;
    Compute $\mathbb{E}_{\mathsf{P}_{d^v}}[Q(y^v,\xi)]$\;
    \If{$\mu^v>\mathbb{E}_{\mathsf{P}_{d^v}}[Q(y^v,\xi)]$ }{
      Identify a violated constraint \eqref{eq:flpref:1st:Qxe} and add it to \eqref{eq:rmp}\;
    }
    \Else{   
      \texttt{solved}$\gets$\texttt{true}\;
    }
    $v\gets v+1$\;
  }
\end{algorithm}

We stress that, in practical implementations, \Cref{alg:ls} is embedded in a Branch \& Cut framework that solves \eqref{eq:rmp}. That is, \Cref{alg:ls} is initiated upon reaching Branch \& Cut nodes which produce an integer feasible solution to \eqref{eq:rmp}. The optimality cuts obtained are, however, globally valid and can hence be broadcast to all open nodes of the Branch \& Cut tree.

\subsection{Tighter Optimality Cut}\label{sec:bilinear_cuts}
In this section, we derive a new optimality cut. 
Throughout, consider a generic iteration $v$ at which the cut is generated. Denote then $d^v$ the distribution induced by the solution to \eqref{eq:rmp} and $(\rho^{sd^vv},\gamma^{sd^vv})_{s\in\mc{S}_{d^v}}$ the dual solutions to $Q(y^v,\xi_s)$, $s\in\mc{S}_{d^v}$.
With these elements we define
\begin{align}
  a^v_{ij}:=\sum_{s\in\mc{S}_{d^v}}\pi_{d^vs}\rho_{ij}^{sd^vv}\xi_{jd^vs}
  \quad\text{ and }\quad
  c_v:=\sum_{s\in\mc{S}_{d^v}}\pi_{d^vs}\sum_{i\in\mc{I}}\gamma_i^{sd^vv}M_i
  \label{eq:cutcoef}
\end{align}
so that the right-hand side of \eqref{eq:oc} (excluding the last term) reads
\begin{align}
  V^v(y):=c_v+\sum_{i\in\mc{I}}\sum_{j\in\mc{J}}a^v_{ij}y_{ij}
  \label{eq:Vv}
\end{align}
Observe that the coefficients \eqref{eq:cutcoef} are nonnegative, because $\rho^{sd^vv}$ and $\gamma^{sd^vv}$ belong to $\mc{R}\subset\mathbb{R}^{\vert\mc{I}\vert\vert\mc{J}\vert}_{\geq 0}\times\mathbb{R}^{\vert\mc{I}\vert}_{\geq 0}$, while probabilities, demands and capacities are nonnegative as well. Observe further that each demand point is assigned to at most one facility.
It follows that $V^v(y)$ ranges in $[c_v,\bar{B}_v]$, where
\begin{align}
  \bar{B}_v:=c_v+\sum_{j\in\mc{J}}\max_{i\in\mc{I}}a^v_{ij}
  \label{eq:Bbar}
\end{align}
We finally let $K_{d}(u):=\vert\mc{J}\vert-\sum_{j\in\mc{J}}u_{j,\zeta(j,d)}$ denote the number of demand points whose allocation is different from the one that enforces $\mathsf{P}_{d}$. Then, \eqref{eq:oc} reads $\mu\leq V^v(y)+K_{d^v}(u)U$.

\begin{proposition}\label{prop:ocv2}
  Let $B_v:=\max\{U,\bar{B}_v\}$. The optimality cut
    \begin{subequations}\label{eq:ocv2}
    \begin{align}
\mu\leq &\sum_{s\in\mc{S}_{d^v}}\pi_{d^vs}\left(\sum_{i\in\mc{I}}\sum_{j\in\mc{J}}\rho_{ij}^{sd^vv}\xi_{jd^vs}y_{ij}+\sum_{i\in\mc{I}}\gamma_i^{sd^vv}M_i\right)+B_v\left(\vert\mc{J}\vert-\sum_{j\in\mc{J}}u_{j,\zeta(j,d^v)}\right)\\
-&\sum_{s\in\mc{S}_{d^v}}\pi_{d^vs}\left(\sum_{i\in\mc{I}}\sum_{j\in\mc{J}}\rho_{ij}^{sd^vv}\xi_{jd^vs}y_{ij}+\sum_{i\in\mc{I}}\gamma_i^{sd^vv}M_i\right)\left(\vert\mc{J}\vert-\sum_{j\in\mc{J}}u_{j,\zeta(j,d^v)}\right)\label{eq:oc:2t}
    \end{align}
    \end{subequations}
    is valid for \eqref{eq:mp}. In addition, it is equivalent to \eqref{eq:oc} when $K_{d^v}(u)=0$, while for $K_{d^v}(u)\geq 1$, it is strictly tighter than \eqref{eq:oc} if $V^v(y)>B_v-U$.\color{black}
\end{proposition}
\begin{proof}
  
  We start by showing validity. Using \eqref{eq:Vv}, cut \eqref{eq:ocv2} reads $\mu\leq V^v(y)+K_{d^v}(u)\left(B_v-V^v(y)\right)$.
    Let $(x,y,u,\mu)$ be feasible for \eqref{eq:mp}. Let $d(u)$ be the distribution enforced by $u$.
    If $d(u)=d^v$ then $K_{d^v}(u)=0$. Then, the right-hand side equals $V^v(y)$.
    Since $(\rho^{sd^vv},\gamma^{sd^vv})\in\mc{R}$ is feasible for the dual of $Q(y,\xi_{d^vs})$ for every $y$, weak duality gives $V^v(y)\geq\mathbb{E}_{\mathsf{P}_{d^v}}Q(y,\xi)\geq\mu$. Here, the last inequality follows from \eqref{eq:flpref:1st:Qx}.
    Othwerwise, if $d(u)\neq d^v$ let $K:=K_{d^v}(u)\geq 1$.
    Because $a^v\geq 0$ and each demand point is assigned to at most one facility, $V^v(y)\leq\bar{B}_v\leq B_v$.
    Rewriting the right-hand side as $B_v+(K-1)\big(B_v-V^v(y)\big)$, we see that both terms are nonnegative increments over $B_v$.
    It follows that
    $$V^v(y)+K\left(B_v-V^v(y)\right)\geq B_v\geq U\geq\mathbb{E}_{\mathsf{P}_{d(u)}}Q(y,\xi)\geq\mu$$
    This confirms that cut \eqref{eq:ocv2} is valid.

    For the comparison with \eqref{eq:oc}, it is easy to see that subtracting the two right-hand sides gives $K_{d^v}(u)\left(B_v-U-V^v(y)\right)$, which is zero when $K_{d^v}(u)=0$ and negative exactly when $V^v(y)>B_v-U$.
\end{proof}

\begin{remark}
  Note that \eqref{eq:oc} and \eqref{eq:ocv2} are both valid but are not ordered.
  It follows that separating both in \Cref{alg:ls} provides, jointly, a tighter cut that either alone.
\end{remark}
\color{black}

Note that the optimality cut \eqref{eq:ocv2} introduces bilinear terms. By expanding the products in \eqref{eq:oc:2t}, we successively obtain
\begin{align}
&\sum_{s\in\mc{S}_{d^v}}\pi_{d^vs}\left(\sum_{i\in\mc{I}}\sum_{j\in\mc{J}}\rho_{ij}^{sd^vv}\xi_{jd^vs}y_{ij}+\sum_{i\in\mc{I}}\gamma_i^{sd^vv}M_i\right)\left(\vert\mc{J}\vert-\sum_{j\in\mc{J}}u_{j,\zeta(j,d^v)}\right)\\
=&\left(\sum_{s\in\mc{S}_{d^v}}\pi_{d^vs}\left(\sum_{i\in\mc{I}}\sum_{j\in\mc{J}}\rho_{ij}^{sd^vv}\xi_{jd^vs}y_{ij}\right)+\sum_{s\in\mc{S}_{d^v}}\pi_{d^vs}\left(\sum_{i\in\mc{I}}\gamma_i^{sd^vv}M_i\right)\right)\left(\vert\mc{J}\vert-\sum_{j\in\mc{J}}u_{j,\zeta(j,d^v)}\right)\\
=&\sum_{s\in\mc{S}_{d^v}}\pi_{d^vs}\left(\sum_{i\in\mc{I}}\sum_{j\in\mc{J}}\rho_{ij}^{sd^vv}\xi_{jd^vs}y_{ij}\right)\left(\vert\mc{J}\vert-\sum_{j\in\mc{J}}u_{j,\zeta(j,d^v)}\right)\\
+&\sum_{s\in\mc{S}_{d^v}}\pi_{d^vs}\left(\sum_{i\in\mc{I}}\gamma_i^{sd^vv}M_i\right)\left(\vert\mc{J}\vert-\sum_{j\in\mc{J}}u_{j,\zeta(j,d^v)}\right)\\
=&\sum_{s\in\mc{S}_{d^v}}\pi_{d^vs}\left(\sum_{i\in\mc{I}}\sum_{j\in\mc{J}}\rho_{ij}^{sd^vv}\xi_{jd^vs}\vert\mc{J}\vert y_{ij}\right)\\
-&\sum_{s\in\mc{S}_{d^v}}\pi_{d^vs}\left(\sum_{i\in\mc{I}}\sum_{j\in\mc{J}}\rho_{ij}^{sd^vv}\xi_{jd^vs} y_{ij}\left(\sum_{j_2\in\mc{J}}u_{j_2,\zeta(j,d^v)}\right)\right)\\
+&\sum_{s\in\mc{S}_{d^v}}\pi_{d^vs}\left(\sum_{i\in\mc{I}}\gamma_i^{sd^vv}M_i\right)\left(\vert\mc{J}\vert-\sum_{j\in\mc{J}}u_{j,\zeta(j,d^v)}\right)\\
=&\sum_{s\in\mc{S}_{d^v}}\pi_{d^vs}\left(\sum_{i\in\mc{I}}\sum_{j\in\mc{J}}\rho_{ij}^{sd^vv}\xi_{jd^vs}\vert\mc{J}\vert y_{ij}\right)\\
-&\sum_{s\in\mc{S}_{d^v}}\pi_{d^vs}\left(\sum_{i\in\mc{I}}\sum_{j\in\mc{J}}\rho_{ij}^{sd^vv}\xi_{jd^vs} \sum_{j_2\in\mc{J}}y_{ij}u_{j_2,\zeta(j,d^v)}\right)\\
+&\sum_{s\in\mc{S}_{d^v}}\pi_{d^vs}\left(\sum_{i\in\mc{I}}\gamma_i^{sd^vv}M_i\right)\left(\vert\mc{J}\vert-\sum_{j\in\mc{J}}u_{j,\zeta(j,d^v)}\right)\\
=&\sum_{s\in\mc{S}_{d^v}}\pi_{d^vs}\left(\sum_{i\in\mc{I}}\sum_{j\in\mc{J}}\rho_{ij}^{sd^vv}\xi_{jd^vs}\vert\mc{J}\vert y_{ij}\right)-\sum_{s\in\mc{S}_{d^v}}\pi_{d^vs}\left(\sum_{i\in\mc{I}}\sum_{j\in\mc{J}}\rho_{ij}^{sd^vv}\xi_{jd^vs} \sum_{j_2\in\mc{J}}y_{ij}u_{j_2,\zeta(j,d^v)}\right)\\
+&\sum_{s\in\mc{S}_{d^v}}\pi_{d^vs}\left(\sum_{i\in\mc{I}}\gamma_i^{sd^vv}M_i \vert\mc{J}\vert\right)-\sum_{s\in\mc{S}_{d^v}}\pi_{d^vs}\left(\sum_{i\in\mc{I}}\gamma_i^{sd^vv}M_i \sum_{j\in\mc{J}}u_{j,\zeta(j,d^v)}\right).
\end{align}
We can now rewrite the optimality cut as follows:
\begin{align}
  \mu\leq &\sum_{s\in\mc{S}_{d^v}}\pi_{d^vs}\left(\sum_{i\in\mc{I}}\sum_{j\in\mc{J}}\rho_{ij}^{sd^vv}\xi_{jd^vs}y_{ij}+\sum_{i\in\mc{I}}\gamma_i^{sd^vv}M_i\right)+B_v\left(\vert\mc{J}\vert-\sum_{j\in\mc{J}}u_{j,\zeta(j,d^v)}\right)\\
  -&\sum_{s\in\mc{S}_{d^v}}\pi_{d^vs}\left(\sum_{i\in\mc{I}}\sum_{j\in\mc{J}}\rho_{ij}^{sd^vv}\xi_{jd^vs}\vert\mc{J}\vert y_{ij}\right)+\sum_{s\in\mc{S}_{d^v}}\pi_{d^vs}\left(\sum_{i\in\mc{I}}\sum_{j\in\mc{J}}\rho_{ij}^{sd^vv}\xi_{jd^vs} \sum_{j_2\in\mc{J}}y_{ij}u_{j_2,\zeta(j,d^v)}\right)\label{eq:oc:bilinearterm}\\
  -&\sum_{s\in\mc{S}_{d^v}}\pi_{d^vs}\left(\sum_{i\in\mc{I}}\gamma_i^{sd^vv}M_i \vert\mc{J}\vert\right)+\sum_{s\in\mc{S}_{d^v}}\pi_{d^vs}\left(\sum_{i\in\mc{I}}\gamma_i^{sd^vv}M_i \sum_{j\in\mc{J}}u_{j,\zeta(j,d^v)}\right)
\end{align}
The bilinear terms have now been isolated in the second term of \eqref{eq:oc:bilinearterm} and involve only products of binary variables. They can be linearized using the Fortet linearization approach \citep{FORTET} described next. 
We replace the products $y_{ij_1}u_{j_2,z}$ by adding to \eqref{eq:rmp} new variables $\ell_{ij_1j_2z}$, for $i\in\mc{I}$, $j_1,j_2\in\mc{J}$, $z\in\mc{Z}$. The new variables are constrained by the Fortet inequalities \eqref{eq:lin_constr}: 
\begin{subequations}\label{eq:lin_constr}
  \begin{align}
    & \ell_{ij_1j_2z}\leq y_{ij_1} & i\in\mc{I}, j_1,j_2\in\mc{J},z\in\mc{Z}\label{eq:lin_constr1}\\
    & \ell_{ij_1j_2z}\leq u_{j_2z} & i\in\mc{I}, j_1,j_2\in\mc{J},z\in\mc{Z}\label{eq:lin_constr2}\\
    & \ell_{ij_1j_2z}+1\geq y_{ij_1}+u_{j_2z} & i\in\mc{I}, j_1,j_2\in\mc{J},z\in\mc{Z} \label{eq:lin_constr3}\\
    & \ell_{ij_1j_2z}+1\geq 0 & i\in\mc{I}, j_1,j_2\in\mc{J},z\in\mc{Z}  \label{eq:lin_constr4}
  \end{align}
\end{subequations}
Constraints \eqref{eq:lin_constr4} are always satisfied when the $\ell$ variables are nonnegative.
The new linear optimality cut can then be reformulated as:
\begin{align}
  \mu\leq &\sum_{s\in\mc{S}_{d^v}}\pi_{d^vs}\left(\sum_{i\in\mc{I}}\sum_{j\in\mc{J}}\rho_{ij}^{sd^vv}\xi_{jd^vs}y_{ij}+\sum_{i\in\mc{I}}\gamma_i^{sd^vv}M_i\right)+B_v\left(\vert\mc{J}\vert-\sum_{j\in\mc{J}}u_{j,\zeta(j,d^v)}\right) \notag\\
  -&\sum_{s\in\mc{S}_{d^v}}\pi_{d^vs}\left(\sum_{i\in\mc{I}}\sum_{j\in\mc{J}}\rho_{ij}^{sd^vv}\xi_{jd^vs}\vert\mc{J}\vert y_{ij}\right) \notag\\ +&\sum_{s\in\mc{S}_{d^v}}\pi_{d^vs}\left(\sum_{i\in\mc{I}}\sum_{j_1\in\mc{J}}\sum_{j_2\in\mc{J}}\rho_{ij_1}^{sd^vv}\xi_{j_1d^vs} \ell_{ij_1j_2\zeta(j_2,d^v)}\right) \notag\\
  -&\sum_{s\in\mc{S}_{d^v}}\pi_{d^vs}\left(\sum_{i\in\mc{I}}\gamma_i^{sd^vv}M_i \vert\mc{J}\vert\right)+\sum_{s\in\mc{S}_{d^v}}\pi_{d^vs}\left(\sum_{i\in\mc{I}}\gamma_i^{sd^vv}M_i \sum_{j\in\mc{J}}u_{j,\zeta(j,d^v)}\right). \label{eq:oc:bilinear}
\end{align}

Observe that this reformulation may require a substantial number of linearization variables and constraints \eqref{eq:lin_constr}. For this reason, we have identified different ways to introduce them. 
The three ways differ in which of the Fortet inequalities \eqref{eq:lin_constr} are imposed on \eqref{eq:rmp} from the outset and which are left to be satisfied otherwise. One point is common to all three and constrains the choice. The variables $\ell_{ij_1j_2z}$ appear in \eqref{eq:oc:bilinear} with a nonnegative coefficient, so the cut is valid only where $\ell_{ij_1j_2z}\geq y_{ij_1}u_{j_2z}$. Constraints \eqref{eq:lin_constr3} are the ones enforcing that inequality, and every version must therefore guarantee them.\color{black}

The first way, which we refer to as \textbf{version 1} (V1), consists of adding the linearization variables $\ell$ to \eqref{eq:rmp} as continuous variables in $[0,1]$. In addition, we add to \eqref{eq:rmp} only the linearization constraints \eqref{eq:lin_constr1} and \eqref{eq:lin_constr2}. The objective function will ensure satisfaction of constraints \eqref{eq:lin_constr3} at the optimal solution. That is, when $y_{ij_1}=u_{j_2z}=1$, the objective will automatically set $\ell_{ij_1j_2z}$ to $1$ since it appears with a nonnegative coefficient in \eqref{eq:oc} (see \Cref{alg:solutionQ} for the dual solution). Because $\ell_{ij_1j_2z}$ is bounded above by both $y_{ij_1}$ and $u_{j_2z}$, this upper bound coincides with the product whenever the latter are integer, so \eqref{eq:lin_constr3} holds at every integer solution and the cut is valid.\color{black}

The second way, which we refer to as \textbf{version 2} (V2), consists of adding the linearization variables $\ell$ upfront to \eqref{eq:rmp}, declared as binary, together with the linearization constraints \eqref{eq:lin_constr3} only. The remaining linearization constraints are added dynamically when found violated. That is, upon reaching integer feasible solutions in the Branch \& Cut that solves \eqref{eq:rmp}, we check whether any of them is violated and we add them as feasibility cuts. To reduce the number of such constraints we separate \eqref{eq:lin_constrl3} in place of \eqref{eq:lin_constr1} and \eqref{eq:lin_constr2}, which is possible because the variables $\ell$ are binary.
\begin{subequations}\label{eq:lin}
  \begin{align}
    & 2\ell_{ij_1j_2z}\leq y_{ij_1}+u_{j_2z} & i\in\mc{I}, j_1,j_2\in\mc{J},z\in\mc{Z}\label{eq:lin_constrl3}.
  \end{align}
\end{subequations}
 \color{black}

Finally, the third way, which we refer to as \textbf{version 3} (V3), consists of adding the linearization variables $\ell$ as binary variables and then adding only \eqref{eq:lin_constrl3} to \eqref{eq:rmp} upfront.
In this case constraints \eqref{eq:lin_constr3} are guaranteed at optimality and constraints \eqref{eq:lin_constr4} are always satisfied. This approach significantly reduces the number of required constraints at the expense of declaring the $\ell$ variables binary. Here \eqref{eq:lin_constrl3} is present at every node, so when we have $y_{ij_1}=u_{j_2z}=1$ the variables $\ell_{ij_1j_2z}$ are free up to one and the objective sets them to one, which is the value of the product. Thus, V3 imposes the upper Fortet direction and obtains the lower one from the objective.

To summarise, V1 imposes \eqref{eq:lin_constr1} and \eqref{eq:lin_constr2} and obtains \eqref{eq:lin_constr3} from the objective. V2 imposes \eqref{eq:lin_constr3} and separates \eqref{eq:lin_constrl3}. V3 imposes \eqref{eq:lin_constrl3} and obtains \eqref{eq:lin_constr3} from the objective. All three carry the same $\vert\mc{I}\vert\vert\mc{J}\vert^2(Z+1)$ variables $\ell$, and, as reported in \Cref{app:size}, comparable numbers of constraints. The main difference between the three versions are the inequalities imposed and those left to be recovered.\color{black}

\subsection{Valid Inequality}\label{sec:validinequality}
 In this section, we propose a valid inequality that is not distribution-specific and speeds up the convergence of the method.

\begin{proposition}\label{prop:vi}
For all $d\in\mc{D}$ and for all $y$, we have
$$\mathbb{E}_{{P}_d}[Q(y,\xi)]\leq \sum_{i\in\mc{I}}\sum_{j\in\mc{J}}R_{ij}\max_{d\in\mc{D}}\min\{M_i,\mathbb{E}_{P_d}[\xi_j]\}y_{ij} \ ,$$
where $\mathsf{P}_d[\xi_j]$ is the marginal expectation of $\xi_j$ according to distribution $d\in\mc{D}$.
\end{proposition}
\begin{proof}
Consider any arbitrary $y$ and $\xi$. Define $w:=(w_{ij})_{i\in\mc{I},j\in\mc{J}}$. The corresponding second-stage problem is:
    \begin{align*}
    Q(y,\xi) = & \max_{w\geq 0} \left\{\sum_{i\in\mc{I}}\sum_{j\in\mc{J}}R_{ij}w_{ij}(\xi)\vert w_{ij}\leq \xi_j y_{ij}, i\in\mc{I},j\in\mc{J},\sum_{j\in\mc{J}}w_{ij}\leq M_i, i\in\mc{I}\right\}.
    \end{align*}
Observe that $Q(y,\xi)$ is concave in $y$ and $\xi$. By Jensen's inequality it follows that, given $y$, for all $d\in\mc{D}$: 
    \begin{align*}
      \mathbb{E}_{{P}_d}[Q(y,\xi)] \leq &    Q(y,\mathbb{E}_{{P}_d}[\xi]) \\
      =& \max_{w\geq 0} \left\{\sum_{i\in\mc{I}}\sum_{j\in\mc{J}}R_{ij}w_{ij}(\xi)\vert w_{ij}\leq \mathbb{E}_{P_d}[\xi_j] y_{ij}, i\in\mc{I},j\in\mc{J},\sum_{j\in\mc{J}}w_{ij}\leq M_i, i\in\mc{I}\right\} \ .
    \end{align*}
    In turn, relaxing the capacity constraints, we obtain:
    \begin{align*}
      &\max_{w\geq 0} \left\{\sum_{i\in\mc{I}}\sum_{j\in\mc{J}}R_{ij}w_{ij}(\xi)\vert w_{ij}\leq \mathbb{E}_{P_d}[\xi_j] y_{ij}, i\in\mc{I},j\in\mc{J},\sum_{j\in\mc{J}}w_{ij}\leq M_i, i\in\mc{I}\right\}\\
      \leq &\max_{w\geq 0} \left\{\sum_{i\in\mc{I}}\sum_{j\in\mc{J}}R_{ij}w_{ij}(\xi)\vert w_{ij}\leq \min\{M_i,\mathbb{E}_{P_d}[\xi_j]\} y_{ij}, i\in\mc{I},j\in\mc{J}\right\}\\
      \leq &\max_{w\geq 0} \left\{\sum_{i\in\mc{I}}\sum_{j\in\mc{J}}R_{ij}w_{ij}(\xi)\vert w_{ij}\leq \max_{d\in\mc{D}}\min\{M_i,\mathbb{E}_{P_d}[\xi_j]\} y_{ij}, i\in\mc{I},j\in\mc{J}\right\}\\
          = &\sum_{i\in\mc{I}}\sum_{j\in\mc{J}}R_{ij} \max_{d\in\mc{D}}\min\{M_i,\mathbb{E}_{P_d}[\xi_j]\} y_{ij} \ .
    \end{align*}
  \end{proof}

By effect of \Cref{prop:vi}, the following inequality is valid for all probability distributions
\begin{align}\label{eq:vi}
  \mu\leq \sum_{i\in\mc{I}}\sum_{j\in\mc{J}}R_{ij} \max_{d\in\mc{D}}\min\{M_i,\mathbb{E}_{P_d}[\xi_j]\} y_{ij}\tag{VI} \, 
\end{align}
and can be inserted to strengthen \eqref{eq:rmp}.

\section{Numerical Experiments}\label{sec:experiments}

We perform numerical experiments on the instances of the problem presented by \citet{LapL94}.
The instances have been extended in such a way to account for decision-dependent uncertain demand.
Following the recipe proposed by \citet{LapL94}, we randomly generate $|\mc{J}|$ demand points with uniform probability in a square of area $[0,100]\times [0,100]$.
Demands $\xi_j$ are assumed to follow truncated normal distributions parametrized by location decisions as shown later.
Each demand point $j\in\mc{J}$ is assigned a \textit{base mean} $\mu_j$ randomly drawn uniformly in the interval $[10,50]$ and \textit{base standard deviation} satisfying a ratio $\sigma_j/\mu_j$ randomly drawn uniformly in $[0.05,0.35]$.
As we will see, base mean and standard deviation will form the basis for decision-dependent distributions.
Three potential facility locations are then generated in the $10\times 10$ square centered around each of the three demand points with the largest $\mu_j$.
The remaining $|\mc{I}|-3$ locations are randomly generated in the $60\times 60$ square centered around the point $(50,50)$.

We assume that the marginal distribution of $\xi_j$ is fully determined by the facility assigned to the demand point $j$.
To mimic a phenomenon by which demand is positively correlated with proximity to the assigned facility, we proceed as follows.
We give subsets $\mc{I}_z$, $z=1,\ldots,Z$ the interpretation of geographical zones.
Geographical zones are created by partitioning facility locations into $Z$ clusters according to their coordinates. We use the \texttt{k-means++} algorithm, see \cite{ArtV07}.

For each demand point, we rank zones $\mc{I}_z$, $z=1,\ldots,Z$ in non-decreasing order of the distance between the demand point and the centroid of the zone.
We let $n(z,j)$ be the rank of zone $z$ in terms of distance from demand point $j$.
Hence, $n(z,j)=1$ indicates that zone $z$ is the closest to demand point $j$, while $n(z,j)=Z$ is the furthest and $n(z,j)=Z+1$ applies to the dummy zone. 
We assume that the demand at point $j$ follows a normal distribution truncated to the left at $0$ with mean and standard deviation depending on the zone to which the demand point is assigned in a given distribution.
Particularly, we set the mean demand for demand point $j$ under distribution $\mathsf{P}_{d}$ to:
$$\mu_{jd}=\mu_j(1+\alpha_{j}^{n(\zeta(j,d),j)})$$
with $\alpha_{j}\in[0,1)$. In all the experiments we set $\alpha_j=0.5$ for every demand point. \color{black} That is, the mean under distribution $\mathsf{P}_{d}$ is equal to the base-level mean $\mu_j$,
increased by $100\alpha_{j}^{n(\zeta(j,d),j)}\%$.
Observe that the exponent ${n(\zeta(j,d),j)}$ becomes larger the furthest zone $\zeta(j,d)$ is from $j$.
Hence, the demand increase becomes exponentially smaller as customer $j$ is assigned to more distant zones.

In a similar way, we assume that the standard deviation of $\xi_j$ under distribution $\mathsf{P}_{d}$ is
$$\sigma_{jd}=\sigma_j(1-\beta_{j}^{n(\zeta(j,d),j)})$$
with $\beta_{j}\in[0,1)$, set to $\beta_j=0.4$ for every demand point in all the experiments.\color{black} \
That is, the standard deviation under distribution $\mathsf{P}_{d}$ is equal to the base-level standard deviation $\sigma_j$ reduced by $100\beta_{j}^{n(\zeta(j,d),j)}\%$.
In this setting, the exponent ${n(\zeta(j,d),j)}$ increases as the assigned zone $\zeta(j,d)$ becomes farther from $j$. Consequently, the reduction in the standard deviation decreases exponentially with the distance between $j$ and its assigned zone.

For each distribution, we generate an empirical distribution by sampling fifty demand realizations (scenarios) of $\xi$ from the underlying truncated normal distributions by inversion.
\begin{remark}
  
  \Cref{prop:ub} and \Cref{prop:vi}  hold for any probability measure, and in particular for the empirical measure induced by the sampled scenarios. The expectations $\mathbb{E}_{\mathsf{P}_d}[\xi_j]$ are accordingly computed as sample means over the scenarios that are actually generated, so that $U$ and \eqref{eq:vi} apply to the specific problem instance that is solved.
\end{remark}
\color{black}

Consistently with \citet{LapL94}, we generate test instances that vary in the number of facilities, demand points and zones as well as in the parameters of the problem.
The instance sizes considered are $(|\mc{I}|,|\mc{J}|,Z)\in\{6,8\}\times \{10,15,20\}\times \{3,4,5\}$.
The parameters are also set according to \citet{LapL94} as follows. We set the net revenue per unit of demand of $j$ satisfied by facility $i$, $R_{ij}=400$ for all $i,j$. We set the capacity of facility $i$ as $M_i=M|\mc{J}|$ for all $i$, with $M=15$.
We set the fixed opening cost for facility $i$ as $F_i=O|\mc{J}|$ for all $i$, with $O=500$. We set the fixed cost for opening link $i-j$ as $C_{ij}=D\Delta_{i,j}$ with $D=270$ and $\Delta_{i,j}$ being the Euclidean distance between $i$ and $j$.
Finally, since instances are randomly generated, we generate ten random instances for every instance size, which yields a total of $10\times 18=180$ instances. \color{black} 
We remark that decision-dependent uncertainty entails working with $(Z+1)^{\vert\mc{J}\vert}$ distributions in each instance.
In our case, this number ranges from $4^{10}$ to $6^{20}$.

\subsection{Results}
\Cref{alg:ls} has been implemented in Java using Cplex 22.1.1 libraries with optimality cuts added only to integer nodes in the Branch \& Bound tree. Subproblems $Q(y,\xi)$ are solved using a specialized procedure proposed by \citet{LapL94} and reported in \Cref{app:solution:Q}, which returns their primal and dual solutions directly, without invoking a linear programming solver. \color{black}
All tests have been run on a machine equipped with two CPUs Intel(R) Xeon(R) Gold 6230 \@ $2.10$GHz. 
\Cref{alg:ls} and all its variants run on a single thread irrespective of the resources made available, and no parallelism is exploited. All tests have been run with Cplex's default settings, including the optimality tolerance of $10^{-4}$.
The only parameter we set was the time limit, which was set to $1800$ seconds in all tests.
The times we report are wall-clock times measured around the call to the solver.
The optimality gap is calculated as $|\mathtt{bestbound}-\mathtt{bestinteger}|/(10^{-10}+|\mathtt{bestinteger}|)$.
\color{black}

The extensive linearized form of the problem has likewise been solved with Cplex 22.1.1\color{black}.
Nevertheless, given the exponential size of the problem the solver could only solve very small instances ($4-5$ customers). Therefore, no result will be presented for the linearized formulation. 

We first analyze the impact of the valid inequality \eqref{eq:vi} by describing its effect on some of the small instances. 
In particular, we used the decomposition algorithm with optimality cuts \eqref{eq:oc}, with and without \eqref{eq:vi}. In \Cref{tab:novivi} we report the following statistics:
\begin{itemize}
\item The average optimality gap (column \texttt{avgGap}) as a percentage. This is calculated as $100*|\mathtt{bestbound}-\mathtt{bestinteger}|/(10^{-10}+|\mathtt{bestinteger}|)$.
\item The percentage of instances solved to a gap smaller than $0.5$\% (column $<0.5\%$).
\item The average solution time (column \texttt{avgTime}). The keyword \texttt{lim} indicates that every run of the corresponding group terminated at the allocated time limit ($1800$ seconds).
\end{itemize}
From \Cref{tab:novivi} it is evident that the valid inequality \eqref{eq:vi} has a crucial impact on the performance of the algorithm.
Without \eqref{eq:vi}, every one of the sixty runs terminates at the allocated time limit without proving optimality.
The corresponding average gaps are between $45.64\%$ and $70.73\%$.
It follows that the bound on $\mu$, which now relies only on the optimality cuts, remains weak and the Branch \& Cut has difficulties pruning. With \eqref{eq:vi} every instance is solved to optimality in less than two seconds on average.\color{black} 

\begin{table}
  \centering
  \caption{Results with and without \eqref{eq:vi} with varying number of facilities, demand points, and zones. }  \label{tab:novivi}
\begin{tabular}{lll|ccc|ccc}
  \toprule
&  &  &  \multicolumn{3}{c}{Without \eqref{eq:vi}}  & \multicolumn{3}{c}{With \eqref{eq:vi}}    \\
  $|\mc{I}|$ & $|\mc{J}|$& $Z$ & \texttt{avgGap}[\%] &   $<0.5\%$   & \texttt{avgTime}[sec.] & \texttt{avgGap}[\%]  & $<0.5\%$ & \texttt{avgTime}[sec.] \\
  \midrule
6 & 10 & 3 & 45.64 & 0\% & \texttt{lim} & 0.00 & 100\% & 0.8\\
6 & 10 & 4 & 57.45 & 0\% & \texttt{lim} & 0.00 & 100\% & 1.0\\
6 & 10 & 5 & 58.98 & 0\% & \texttt{lim} & 0.00 & 100\% & 2.0\\
8 & 10 & 3 & 53.90 & 0\% & \texttt{lim} & 0.00 & 100\% & 0.6\\
8 & 10 & 4 & 59.56 & 0\% & \texttt{lim} & 0.00 & 100\% & 1.1\\
8 & 10 & 5 & 70.73 & 0\% & \texttt{lim} & 0.00 & 100\% & 1.9\\
  \bottomrule
\end{tabular}
\end{table}

Part of the explanation behind the impact of \eqref{eq:vi} is as follows.
In this problem, the space of feasible solutions $\left\{(x_i)_{i\in\mc{I}}\in\{0,1\}^{\vert\mc{I}\vert},(y_{ij})_{i\in\mc{I},j\in\mc{J}}\in\{0,1\}^{\vert\mc{I}\vert\times \vert\mc{J}\vert}\right\}$ is essentially partitioned into $\vert\mc{D}\vert$ subsets.
Each subset contains the solutions $(x,y)$ that correspond to a specific allocation of demand points to subsets of facilities. Thus, each subset generates a different distribution of the demand.
The optimality cuts \eqref{eq:oc} reconstruct the expected second-stage revenue function in each of these subsets.
The number of subsets is, however, exponentially large in $\vert\mc{J}\vert$ (recall that in \Cref{tab:novivi} we use $\vert\mc{J}\vert=10$).
Therefore, the algorithm needs to generate a number of cuts sufficient to  reconstruct a very large number of expectations. The main effect of \eqref{eq:vi} is that it provides non-trivial decision-independent bounds that are valid across multiple subsets of the feasible space.
These, in turn, simplify bounding and pruning in the Branch \& Cut algorithm. For this reason, all the results presented in what follows are obtained using \eqref{eq:vi}.

We now compare DD-LS, that is the algorithm with the original optimality cuts \eqref{eq:oc}, with the three versions V1--V3 of the algorithm that, additionally, separate the new cut \eqref{eq:oc:bilinear}.
Recall from \Cref{sec:bilinear_cuts} that the two inequalities are both valid and that neither dominates the other, so that the versions V1--V3 separate both at every round. The comparison is therefore between the original cut alone and the combination of the two.
The results are provided in \Cref{tab:r2:main}.
They are grouped by number of demand points, which is the main driver both of the number of distributions, $(Z+1)^{\vert\mc{J}\vert}$, and of the number of linearization variables $\ell_{ij_1j_2z}$.
Each cell of \Cref{tab:r2:main} aggregates the $60$ runs obtained by combining $\vert\mc{I}\vert\in\{6,8\}$, $Z\in\{3,4,5\}$ and the ten instances generated for each size.

\begin{table}[]
    
    \centering
    \caption{Comparison of DD-LS with versions V1--V3 of the algorithm using the new optimality cut \eqref{eq:oc:bilinear} together with \eqref{eq:oc}, by number of demand points. Each cell aggregates the $60$ runs obtained by combining $\vert\mc{I}\vert\in\{6,8\}$, $Z\in\{3,4,5\}$ and the ten instances generated for each size.}\label{tab:r2:main}
    \begin{tabular}{llccc}
\toprule
 & $\vert\mc{J}\vert$ & \texttt{solved} & \texttt{avgTime}[sec] & \texttt{avgGap}[\%] \\
\midrule
\multirow[c]{3}{*}{DD-LS} & 10 & 52/60 & 1.3 & 0.00 \\
 & 15 & 31/60 & 12.2 & 0.00 \\
 & 20 & 29/60 & 167.9 & 0.27 \\
\midrule
\multirow[c]{3}{*}{V1} & 10 & 52/60 & 6.3 & 0.00 \\
 & 15 & 37/60 & 146.5 & 0.12 \\
 & 20 & 25/60 & 501.9 & 1.01 \\
\midrule
\multirow[c]{3}{*}{V2} & 10 & 52/60 & 11.9 & 0.00 \\
 & 15 & 39/60 & 158.9 & 0.09 \\
 & 20 & 27/60 & 544.4 & 1.20 \\
\midrule
\multirow[c]{3}{*}{V3} & 10 & 48/60 & 74.8 & 0.12 \\
 & 15 & 19/60 & 396.9 & 0.74 \\
 & 20 & 20/60 & 891.1 & 2.91 \\
\bottomrule
\end{tabular}
\end{table}

On the instances with $15$ demand points the combination is clearly stronger than the original cut alone: V2 closes $39$ instances out of $60$ and V1 closes $37$, against $31$ for DD-LS. Over the whole test bed V2 solves $118$ instances out of $180$ and V1 solves $114$, against $112$ for DD-LS. The benefit is, however, narrow.
DD-LS is faster at every size, and by a wide margin.
It also leaves smaller gaps, $0.27\%$ against $1.01\%$, $1.20\%$ and $2.91\%$ on the instances with $20$ demand points, where it additionally closes more instances than any of the three versions. V3, which imposes the largest number of linearization constraints upfront, is the weakest of the three throughout.
A practitioner with a fixed time budget would therefore choose the original cut, and would turn to the combination only when the objective is to close the gap on instances of moderate size that the original cut leaves open.

An explanation for the gap left by V3 on the largest instances, and for the poor performance of V1 and V2, is to be found in the size of \eqref{eq:rmp} once the linearization variables and constraints are added.
That size is reported in \Cref{app:size}.
The central issue is the number of bilinear terms $y_{ij_1}u_{j_2z}$, each of which is replaced by a variable $\ell_{ij_1j_2z}$.
Their number is grows quadratically in the number of demand points.
It affects all three versions since each of them adds the same $\ell$ variables.
With $\vert\mc{J}\vert=20$ and $Z=5$ \eqref{eq:rmp} carries on average $17068$ variables against the $268$ of DD-LS.
V1 imposes two inequalities per product and reaches $33881$ constraints on average, while V2 and V3 impose one and reach $17082$ and $17081$ respectively, against $281$ for DD-LS.
What distinguishes the three versions in practice is therefore not the size of \eqref{eq:rmp}, which is comparable, but which of the Fortet inequalities is imposed and which is separated, and the consequences this has for how quickly the search can make progress.
\color{black}

\begin{table}[]
    
    \centering
    \caption{Quality of the bounds.
      Each row is a separate pairwise comparison between DD-LS and one version, restricted to the \texttt{n} instances that neither of the two solved to optimality.
      Hence, the three rows refer to three different sets of instances.
      In each pair of columns, DD-LS is the algorithm with \eqref{eq:oc} and V* is the version named in the row, that is V1, V2 or V3.
      \texttt{nCutRounds} is the number of cut separation rounds, each adding one cut for DD-LS and two for V* (i.e., cut \eqref{eq:oc:bilinear} together with \eqref{eq:oc}).
      \texttt{rootDist} and \texttt{endDist} are the relative distances of the bound at the end of the root node and at termination, respectively, from the best solution found by any of the methods.
      That is, both are averages, over the \texttt{n} instances of the row, of $100(z^*-L)/\vert z^*\vert$, where $z^*$ is the value of the best solution found on that instance by any of the methods and $L$ is the bound in question.
      Note that $z^*$ is common to the two methods being compared, so that only their bounds differ. Furthermore, $z^*$ is not the incumbent of either of them.}\label{tab:r2:bounds}
    \begin{tabular}{lcccccccc}
\toprule
 & & \multicolumn{2}{c}{\texttt{nCutRounds}} & \multicolumn{2}{c}{\texttt{rootDist}[\%]} & \multicolumn{2}{c}{\texttt{endDist}[\%]} \\
\cmidrule(lr){3-4}\cmidrule(lr){5-6}\cmidrule(lr){7-8}
 & \texttt{n} & DD-LS & V* & DD-LS & V* & DD-LS & V* \\
\midrule
V1 & 53 & 1581 & 555 & 4.64 & 4.83 & 0.31 & 1.11 \\
V2 & 49 & 1698 & 554 & 4.73 & 4.75 & 0.34 & 1.23 \\
V3 & 62 & 1360 & 187 & 4.34 & 4.41 & 0.27 & 2.04 \\
\bottomrule
\end{tabular}
\end{table}

To better explain the results it is useful to assess the difference in terms of dual bounds, which we report in \Cref{tab:r2:bounds}.
The comparison made there is pairwise in the sense that each row of the table compares DD-LS with one single version, and is therefore restricted to the instances that neither of those two closes to optimality.
On those instances the bound at the end of the root node is essentially the same for the two methods being compared.
However, the bound at termination is significantly better for DD-LS.
The reason is visible in the \texttt{nCutRounds} columns.
Each round of separation is more expensive once the bilinear terms are present and the linear programming relaxations involved are larger due to the linearization variables and constraints.
It follows that, within the same time limit, DD-LS performs roughly three times as many separation rounds as V1 and V2, and more than seven times as many as V3.
An individual combined cut applied by V1-V3 is tighter, in the sense made precise in \Cref{prop:ocv2}, but the algorithm is able to generate significantly fewer of them.
It follows that the reconstruction of the expected recourse function proceeds more slowly. 

Finally, the results illustrate that separating \eqref{eq:oc} alongside \eqref{eq:oc:bilinear}, rather than the latter alone, is preferable.
When we drop \eqref{eq:oc}, V2 solves $43$ fewer instances, and more than triples its average solution time.
V3 drops from $87$ to $59$ and V1 from $114$ to $105$ solved instances.
Since \eqref{eq:oc} involves none of the linearization variables, it is inexpensive to separate, and it bounds $\mu$ on exactly those solutions for which the bilinear cut is weak.
\color{black}

\section{Conclusions}\label{sec:conclusions}
In this article, we introduce a model that considers decision-dependent uncertain demand in single-source capacitated facility location problems.
The model is built on the idea of having the distribution of demand dependent on the facility to which demand points are allocated.
By partitioning the set of facilities, the model can capture, at one extreme, the case where demands are not influenced by allocation decisions, and, at the other extreme, the case where each individual allocation decision determines a different distribution.
Between these extremes, one can find cases where the distribution depends on the subset of facilities to which the demand point is allocated.

To address the exponential size of the model and the possible number of probability distributions, we propose a cutting-plane algorithm based on the L-Shaped method for problems with decision-dependent uncertainty (see \cite{PanH25}), which generates distribution-specific cuts. 
In addition, we derive a new, tighter, optimality cut based on a bilinear reformulation of the original cut as well as a tight valid inequality.

Extensive numerical tests demonstrate that the method, supported by the valid inequality, can find optimal or near optimal solutions to problem instances with up to $6^{20}$ potential distributions.
This improves significantly the current ability to solve location problems with decision-dependent uncertainty.
We also assess the practical value of the new optimality cut.
Separated alongside the original one it closes slightly more instances.
However, the original cut remains preferable overall, both in terms of solution time and optimality gap attained.
The main reason is that the bilinear terms make each round of separation considerably more expensive and the algorithm therefore affords substantially fewer cuts within the same time limit.\color{black}

\section{Acknowledgments}
This study was partially funded by the Novo Nordisk Foundation grant NNF24OC0089770.

\bibliographystyle{plainnat} 
\bibliography{preprint.bib}

\begin{appendices}
    \crefalias{section}{appsec}
\end{appendices}
\appendix
\section{Size of the Relaxed Master Problem}\label{app:size}

\begin{table}
    \centering
    \caption{Size of Relaxed Master Problem for DD-LS and Versions V1-V3 of Algorithm with New Optimality Cuts \eqref{eq:oc:bilinear}, averaged over the ten instances of each size. The figures refer to the model as built, before any cut is separated, and are therefore the same whether or not \eqref{eq:oc} is added alongside \eqref{eq:oc:bilinear}. Note that V2 now posts \eqref{eq:lin_constr3} upfront, for the reason given in \Cref{sec:bilinear_cuts}, and so no longer starts from a master problem of the size of the one of DD-LS.\color{black}}
    \label{tab:mp_size}
    \begin{tabular}{lllcccccc}
\toprule
 &  &  & \multicolumn{2}{c}{\texttt{\#\_variables}} & \multicolumn{2}{c}{\texttt{\#\_binary}} & \multicolumn{2}{c}{\texttt{\#\_constraints}} \\
 & $\vert\mc{Z}\vert$ & $\vert\mc{J}\vert$& \texttt{avg} & \texttt{max} & \texttt{avg} & \texttt{max} & \texttt{avg} & \texttt{max}  \\
\midrule
\multirow[c]{9}{*}{DD-LS} & \multirow[c]{3}{*}{3} & 10 & 118.00 & 129 & 117.00 & 128 & 121.00 & 131 \\
 &  & 15 & 173.00 & 189 & 172.00 & 188 & 181.00 & 196 \\
 &  & 20 & 228.00 & 249 & 227.00 & 248 & 241.00 & 261 \\
\cline{2-9}
 & \multirow[c]{3}{*}{4} & 10 & 128.00 & 139 & 127.00 & 138 & 131.00 & 141 \\
 &  & 15 & 188.00 & 204 & 187.00 & 203 & 196.00 & 211 \\
 &  & 20 & 248.00 & 269 & 247.00 & 268 & 261.00 & 281 \\
\cline{2-9}
 & \multirow[c]{3}{*}{5} & 10 & 138.00 & 149 & 137.00 & 148 & 141.00 & 151 \\
 &  & 15 & 203.00 & 219 & 202.00 & 218 & 211.00 & 226 \\
 &  & 20 & 268.00 & 289 & 267.00 & 288 & 281.00 & 301 \\
\midrule
\multirow[c]{9}{*}{V1} & \multirow[c]{3}{*}{3} & 10 & 2918.00 & 3329 & 117.00 & 128 & 5721.00 & 6531 \\
 &  & 15 & 6473.00 & 7389 & 172.00 & 188 & 12781.00 & 14596 \\
 &  & 20 & 11428.00 & 13049 & 227.00 & 248 & 22641.00 & 25861 \\
\cline{2-9}
 & \multirow[c]{3}{*}{4} & 10 & 3628.00 & 4139 & 127.00 & 138 & 7131.00 & 8141 \\
 &  & 15 & 8063.00 & 9204 & 187.00 & 203 & 15946.00 & 18211 \\
 &  & 20 & 14248.00 & 16269 & 247.00 & 268 & 28261.00 & 32281 \\
\cline{2-9}
 & \multirow[c]{3}{*}{5} & 10 & 4338.00 & 4949 & 137.00 & 148 & 8541.00 & 9751 \\
 &  & 15 & 9653.00 & 11019 & 202.00 & 218 & 19111.00 & 21826 \\
 &  & 20 & 17068.00 & 19489 & 267.00 & 288 & 33881.00 & 38701 \\
\midrule
\multirow[c]{9}{*}{V2} & \multirow[c]{3}{*}{3} & 10 & 2918.00 & 3329 & 2917.00 & 3328 & 2922.00 & 3332 \\
 &  & 15 & 6473.00 & 7389 & 6472.00 & 7388 & 6482.00 & 7397 \\
 &  & 20 & 11428.00 & 13049 & 11427.00 & 13048 & 11442.00 & 13062 \\
\cline{2-9}
 & \multirow[c]{3}{*}{4} & 10 & 3628.00 & 4139 & 3627.00 & 4138 & 3632.00 & 4142 \\
 &  & 15 & 8063.00 & 9204 & 8062.00 & 9203 & 8072.00 & 9212 \\
 &  & 20 & 14248.00 & 16269 & 14247.00 & 16268 & 14262.00 & 16282 \\
\cline{2-9}
 & \multirow[c]{3}{*}{5} & 10 & 4338.00 & 4949 & 4337.00 & 4948 & 4342.00 & 4952 \\
 &  & 15 & 9653.00 & 11019 & 9652.00 & 11018 & 9662.00 & 11027 \\
 &  & 20 & 17068.00 & 19489 & 17067.00 & 19488 & 17082.00 & 19502 \\
\midrule
\multirow[c]{9}{*}{V3} & \multirow[c]{3}{*}{3} & 10 & 2918.00 & 3329 & 2917.00 & 3328 & 2921.00 & 3331 \\
 &  & 15 & 6473.00 & 7389 & 6472.00 & 7388 & 6481.00 & 7396 \\
 &  & 20 & 11428.00 & 13049 & 11427.00 & 13048 & 11441.00 & 13061 \\
\cline{2-9}
 & \multirow[c]{3}{*}{4} & 10 & 3628.00 & 4139 & 3627.00 & 4138 & 3631.00 & 4141 \\
 &  & 15 & 8063.00 & 9204 & 8062.00 & 9203 & 8071.00 & 9211 \\
 &  & 20 & 14248.00 & 16269 & 14247.00 & 16268 & 14261.00 & 16281 \\
\cline{2-9}
 & \multirow[c]{3}{*}{5} & 10 & 4338.00 & 4949 & 4337.00 & 4948 & 4341.00 & 4951 \\
 &  & 15 & 9653.00 & 11019 & 9652.00 & 11018 & 9661.00 & 11026 \\
 &  & 20 & 17068.00 & 19489 & 17067.00 & 19488 & 17081.00 & 19501 \\
\bottomrule
\end{tabular}
\end{table}

\section{Solution of $Q(y,\xi)$}\label{sec:efficiency}\label{app:solution:Q}
The primal and dual solutions to $Q(y,\xi)$ can be efficiently obtained using a direct method devised by \citet{LapL94}.
The method is based on a specialized sorting procedure and is described in \Cref{alg:solutionQ}.

\begin{algorithm}
\caption{Procedure to find primal and dual solutions to $Q(y,\xi)$.}
  \label{alg:solutionQ}
  \SetKwInOut{Input}{input}\SetKwInOut{Output}{output}
  \Input{A solution $(x^v,y^v)$ to \eqref{eq:rmp} and the corresponding distribution $d^v$}
  \Output{Primal solution $w_{ijs}^v$ and dual multipliers $\rho^{sd^vv}_{ij}$ and $\gamma_i^{sd^vv}$}
  \For{$s\in\mc{S}_{d^v}$}{
    \For{$i\in\mc{I}:x_i^v=1$}{
      Rank demand points $j\in\mc{J}$ in non-increasing order of the values $R_{ij}$\;
      Assign demands $\xi_{jd^vs}$ to $i$ (populating variables $w_{ijs}^v$) in the same order as the values $R_{ij}$
      as long as $y_{ij}^v=1$, $R_{ij}>0$ and the residual capacity $M_i-\sum_{j\in\mc{J}}w_{ijs}^v$ is not exceeded\;
      Let $\lambda_{is}$ be the index $j\in\mc{J}$ such that $y_{ij}^v=1$, $R_{ij}>0$ and the demand of $j$ can be satisfied only partially by $i$.
      If no such index exists, let $\lambda_{is}=|\mc{J}|+1$ with $R_{i,\lambda_{is}}=0$\;
      Let $\rho_{ij}^{sd^vv}=R_{ij}-R_{i,\lambda_{is}}$ if $j<\lambda_{is}$ (on the ranked list), $0$ otherwise\;
      Let $\gamma_i^{sd^vv} = R_{i,\lambda_{is}}$\;
    }
  }
  
\end{algorithm}

\end{document}